УДК 378:51(09)

**БОГАТОВ ЕГОР МИХАЙЛОВИЧ**
**КОРЕНЕВ АРТЕМ ВИКТОРОВИЧ**
**МИХАЙЛОВ ИЛЬЯ СЕРГЕЕВИЧ**

Россия, г. Старый Оскол, СТИ НИТУ МИСиС; г. Губкин, ГФ НИТУ МИСиС

email: embogatov@inbox.ru

**О СОВРЕМЕННЫХ ИНСТРУМЕНТАХ И МЕТОДАХ ВЕДЕНИЯ НАУЧНЫХ ИССЛЕДОВАНИЙ ПО ИСТОРИИ МАТЕМАТИКИ[1]**

*Предложен один из вариантов систематизации деятельности историка математики, а также схема организации исследовательской и поисковой работы при подготовке научных статей и докладов по истории науки.*

*Ключевые слова: схема организации исследований по истории математики, методика подготовки статей и докладов по истории математики, современный инструментарий историко-математических исследований; источниковая база историка математики.*

> ...История наук даёт лучший и наиболее
> надёжный материал, на котором могут
> быть изучены закономерности
> в развитии человечества.
> *В. Оствальд*

**Введение**

До конца XX века историки науки (как отечественные, так и зарубежные) располагали лишь библиотечным фондом на бумажном носителе. В России это: Российская государственная библиотека, библиотеки ведущих университетов, архивы РАН и т.д. Несмотря на это, интерес к истории физико-математических наук был довольно высоким: издавались серьезные журналы: ИМИ (Историко-математические исследования); ИАИ (Историко-астрономические исследования); ИМЕН (История и методология естественных наук); труды ИИЕТ РАН[2] (Институт истории естествознания и техники имени С. И.

---

[1] Работа является расширенным дополненным вариантом доклада [1].
[2] Значительная часть этих изданий доступна на сайте В.Е. Пыркова [2].

Вавилова РАН); выходили энциклопедии по истории математики [3]-[5] и механики [6]. Наиболее известным и полным (на конец 1960-х гг.) изданием по праву считается «История отечественной математики» [7]- [10].

К концу XX - началу XXI в. в России история математики стала понемногу сдавать свои позиции; долгое время она отсутствовала даже в рубрикаторе РФФИ (ситуация поменялась только в конце 2010-х гг.). Это привело к нарушению непрерывности передачи опыта организации историко-математических исследований от старшего поколения к младшему, особенно в научных коллективах, расположенных на периферии. С другой стороны, с развитием сети Интернет появились новые возможности для доступа к первоисточникам и статьям по истории науки, опубликованным на иностранных языках.

Отметим, что история математики помогает решить проблемы мотивации, особенно это касается нынешнего поколения. Если молодые исследователи будут равняться на кого-то из великих людей, внесших вклад в развитие науки (в частности математики), то у них будет дополнительный стимул к получению высоких результатов в своей области.

Кроме того, работы по истории математики всегда имеют пропагандистский характер. Они в большей или меньшей степени демонстрируют важность математики, и как составной части общей культуры, и как инструмента для решения прикладных задач

История математики способствует развитию таких навыков, как коммуникация и критическое мышление, что может стимулировать преподавание модулей истории математики в ВУЗах. В некоторых странах она активно используется для мотивации изучения математики школьниками и студентами [11].

Настоящая работа имеет своей целью систематизировать современный инструментарий историко-математических исследований, а также сориентировать начинающего исследователя в области историко-математических журналов; конференций и литературы справочного характера. Будут предложены рекомендации молодым учёным по написанию статей по истории естественных наук.

## 1.Схема организации исследований

После того, как обозначена постановка задачи (обычно это делает научный руководитель), целесообразно сделать шаги в разных направлениях:

*a) Детальное изучение предметной области.*

Постараться разобраться в предмете как можно лучше. Здесь на помощь может прийти хороший учебник, брошюра из серии «Популярные лекции по математике» издательства «Наука», научно-популярные статьи в журнале «Квант» и в сборнике «Математическое просвещение» [12], материалы летних/зимних математических школ, в том числе, представленных в виде лекций на канале YouTube. Например, если нас интересует история выпуклости, полезно прочесть статью В.М. Тихомирова *Геометрия выпуклости* [13].

Большую ценность для историка науки представляет собой общероссийский портал Math-Net.Ru. Он выполнен в виде современной информационной системы, предоставляющей российским и зарубежным ученым различные возможности в поиске научной информации по математике и другим наукам. Кроме печатных источников, портал предлагает большое количество разнообразных видеороликов по математике, расположенных во вкладке «Видеотека» (например, по теореме о неподвижной точке можно рекомендовать лекцию [14]).

Значимая информация может быть почерпнута также из обзоров последних десятилетий в той или иной области математики, опубликованных в журнале УМН (см. например обзор "О некоторых основных направлениях общей топологии" [15]).

*б) Поиск материалов на схожую тему исследования.*

Посмотреть выпуски ИМИ [16] на предмет наличия статей на близкую к нему тему, а также осуществить поиск в Интернете. Например, сформулировать запрос с помощью Google как на русском, так и на английском языках. В случае с историей выпуклости, к примеру, мы набираем history of convexity и среди ответов видим большую обзорную статью П.М. Грубера в учебнике [17], в которой много ссылок. После этого имеет смысл пристально изучить найденные источники.

Если нам известен хотя бы один из математиков, который занимался интересующей нас темой (например, Архимед), то мы можем воспользоваться интернет-архивом по истории математики MacTutor [18], созданным учеными из университета Сент-Эндрюс в Шотландии Джоном О'Коннором и Эдмундом Робертсоном. Во вкладке BIOGRAPHIES мы обнаружим большую обзорную статью о результатах Архимеда со ссылками на работы историков науки, в том числе отечественных.

На этом этапе можно воспользоваться статьями к юбилеям [19], научными автобиографиями [20] и некрологами [21].

Неоспоримую помощь в поиске необходимых материалов сможет оказать старейший научный архив России - архив Российской академии наук [22]. Он представляет совокупность документов, образующихся в деятельности Российской академии наук и её учреждений, имеющих научное, социально-культурное и историческое значение. В этом архиве, например, имеются письма известных ученых к своим коллегам, из которых иногда бывает можно почерпнуть мотивацию к возникновению той или иной идеи или метода.

Параллельно с работой над статьями по истории математики (п.б)) можно обратиться к энциклопедиям и хрестоматиям таким, как «The Oxford Handbook of The History of Mathematics» [23], «Companion Encyclopedia of the History and Philosophy of the Mathematical Sciences» [24], «История математики от Декарта до середины XIX столетия» [25], «Краткий очерк истории математики» [26], а также к книгам по истории отдельных разделов математики таких, как история топологии [27].

Следует отдельно отметить 3-х томный сборник «Математика, ее содержание, методы и значение» [28-30], изданный в 1956 году в СССР, который знакомит читателя с с содержанием и методами отдельных математических дисциплин, их основами и путями развития.

Для изучения вклада зарубежных ученых помощью будет обращение к путеводителю по историко-математической литературе на английском языке [31].

*в) Использование форумов и научных социальных сетей для поиска ответов на возникшие вопросы.*

В сети Интернет достаточно большое количество форумов, на которых обсуждаются вопросы по истории математики. Можно выделить несколько из них:

- Math10.com.

   Данный русскоязычный форум носит любительский характер. На нем обсуждается возникновение математики, предпосылки ее зарождения, а также последовательность ее развития. Поднимаются, в частности, вопросы о первенстве возникновения «числа» или письменности URL: https://www.math10.com/ru/forum/viewforum.php?f=40

- Readit. History of Mathematics?

   Этот форум представлен на английском языке. На нем выносятся на обсуждение отдельные моменты развития математики, важные и интересные этапы ее становления. Также на нем можно найти полезные ссылки на книги по истории математики.
   URL: https://www.reddit.com/r/history/comments/a16xoj/history_of_mathematics/

- MATHEMATICS.

   Этот форум также представлен на английском языке. На нем поднимаются самые разнообразные вопросы, посвященные истории математики (от причин обозначения той или иной величины определенной буквой до поиска книг по истории математики).
   URL: https://math.stackexchange.com/questions/tagged/math-history

Помимо форумов, существуют специальные социальные сети для общения ученых различных стран. В качестве примера можно привести ResearchGate – европейскую некоммерческую социальную сеть, в рамках которой исследователи могут делиться статьями, задавать вопросы и отвечать на них и находить соавторов для научной работы. Эта сеть предоставляет такие приложения, как семантический поиск (поиск по аннотации), совместное использование файлов, обмен базой публикаций, форумы, методологические дискуссии и так далее.

    г) *Поиск и изучение первоисточников*

После изучения предмета исследований и истории вопроса а)-в) весьма важно обратиться к первоисточникам. Статьи на русском языке можно искать на общероссийском математическом портале Math-Net.ru. Несомненную пользу принесёт обращение к электронной библиотеке «Научное наследие России» [32]. Проект электронной библиотеки "Научное наследие России" нацелен на обеспечение решения задачи о сохранении научного наследия и создании условий его эффективного освоения.

Эта библиотека создавалась учреждениями РАН как общедоступная с целью предоставить пользователям Интернет информацию о выдающихся российских ученых, внесших вклад в развитие фундаментальных естественных наук, и полных текстов опубликованных ими наиболее значительных работ. На этом сайте, в частности, можно найти такие редкие издания, как сборник трудов П.И. Чебышева.

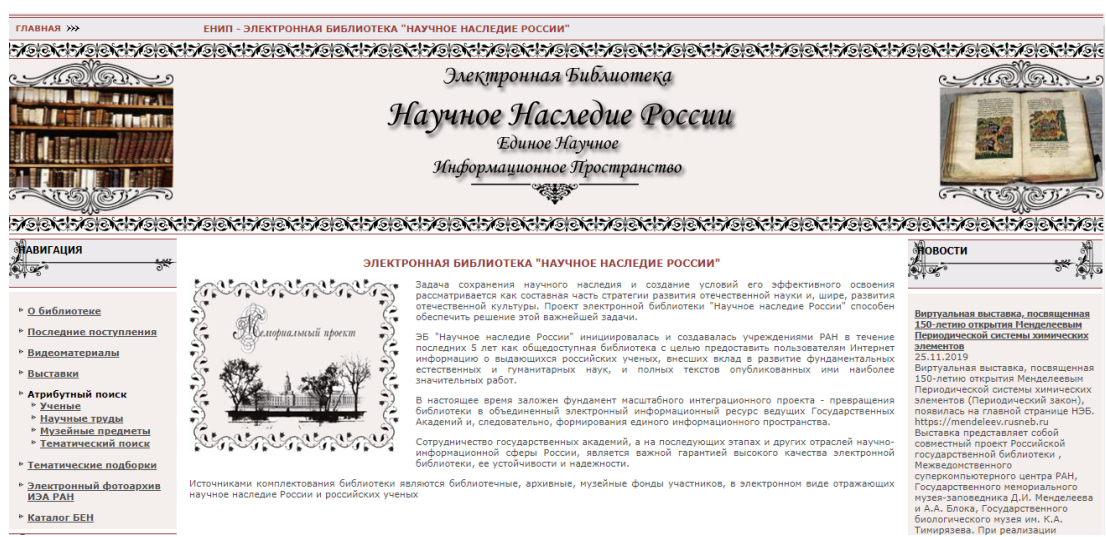

Рис.1 Внешний вид сайта электронной библиотеки "Научное наследие России"

Помимо этого, в поиске первоисточников может помочь электронный ресурс *Project Euclid* [33]. Первоначально он был создан для того, чтобы обеспечить платформу для небольших научных издателей журналов по математике и статистике, а также чтобы перейти от печатной к электронной форме экономически эффективным способом. Благодаря сочетанию поддержки со стороны подписавшихся библиотек и участвующих издателей более 70% журнальных статей, размещенных на сайте Project Euclid, находятся в открытом доступе. В настоящее время целью Project Euclid является предоставление издательских услуг для теоретической и прикладной математики и статистики по всему миру. В рамках данного проекта возможно опубликование препринтов - предпечатных версий статей с целью обозначить приоритет в определённом направлении исследований.

Статьи на других (европейских) языках можно найти на сайтах:
http://www.digizeitschriften.de (немецкая цифровая библиотека Геттингена)
Предоставляется доступ к более чем 700 известным немецким академическим журналам из 21 предметной области с 9 миллионами страниц. Все статьи доступны для скачивания в формате PDF и печати.
https://gallica.bnf.fr (французская национальная библиотека)
Французский архив ставит своей целью оцифровать все содержимое Национальной библиотеки Франции — более 12 миллионов книг и манускриптов, 500 тысяч периодических изданий. В рамках проекта выложено в интернет около 80 тысяч работ и 70 тысяч изображений, идёт работа по оцифровке архива французских статей XIX в. В интерфейсе есть поддержка русского языка.
https://eudml.org (европейская цифровая математическая библиотека)
Цель проекта The European Digital Mathematics Libriaryu (EuDML) — создать единую информационную систему математической литературы, опубликованной в Европейских изданиях. Проект осуществляется при поддержке ведущих Европейских университетов и научных организаций. EuDML индексирует более 250 тысяч публикаций из 14 коллекций журналов. В интерфейсе есть поддержка русского языка
https://www.jstor.org (международная цифровая библиотека академического контента)
Journal STORage предоставляет доступ к более чем 12 миллионам научных журнальных статей, книг и первоисточников по 75 дисциплинам. JSTOR - мощная исследовательская и преподавательская платформа, с помощью которой можно изучить широкий спектр научного контента, в том числе и по истории математики.
https://archive.org (американский цифровой архив)
Этот архив содержит редкие и ценные издания XIX-XXв (книги и статьи), опубликованные, в основном на английском языке. Доступ к некоторым историко-математическим материалам XXI века (таким, как, например, монография Дж. Лютцена о Лиувилле [34]) осуществляется по подписке.

Для того чтобы приступить к переводу статьи на русский язык машинным способом (например, с помощью переводчиков Google или Deepl), находящимися в свободном доступе, необходимо преобразовать файлы pdf в текст. Это можно сделать с помощью интернет-сервисов распознания текста (http://tools.pdf.24.org; https://www.onlineocr.net/ и т.п.)

*д) Воссоздание эволюции математических идей*

Самое трудное – это воссоздание эволюции той или иной математической идеи. В отличие от научного обзора, упор здесь должен делаться не столько на то, *как* появилась и развилась та или иная идея, сколько на то, *почему* ее развитие происходило тем или иным образом. Очень важно видеть историю математики в целом, как большой и связный процесс, и если мы выделим какую-то линию в этом процессе, то необходимо зафиксировать её связи и взаимодействие с другими линиями и процессами.

При работе в области истории математики полезно ставить перед собой вопросы следующего характера.[3]

«Почему возникла та или иная область математики? Какие задачи решались с её помощью? Какую роль при этом играли конкретные разделы обсуждаемой области? Что из обсуждаемых идей пришло из практики, а что решало внутри математические задачи? Как это конфликтовало между собой? Что стояло за теми или иными поворотами в данном научном направлении?»

Суть методологии историко-научных исследований очень образно представил академик Б.В. Гнеденко в докладе на секции истории математики IV Всесоюзного математического съезда: «…*для правильности исторической перспективы надо любое явление в истории науки рассматривать не только в его генезисе, с точки зрения предшествующего, но и с точки зрения последующего, т.е. с учетом его влияния и его последствий, которые следует прослеживать до наших дней. Чтобы охарактеризовать плод, надо знать и дерево, на котором он созрел, и семена, которые он в себе содержит, и ростки, которые дали эти семена*» [35, с. 20].

*е) Апробация*

Предположим, что первая версия статьи по истории математики создана. Естественно, она не будет являться окончательным вариантом. Для апробации таких работ и их улучшения в учёном сообществе существует неписаное правило: *научный результат должен пройти*

---

[3] Вопросы воспроизводятся из переписки доцента СТИ НИТУ МИСиС Е.М. Богатова и профессора МГУ А.В. Боровских, электронное письмо от 26.03.2019.

*обсуждение.* Это обсуждение происходит либо на семинарах, либо на конференциях. По истории математики в России известны два постоянно работающих семинара:

- в МГУ при кабинете истории и методологии математики под руководством президента Международной академии истории науки профессора С.С. Демидова;
- в ПОМИ РАН под руководством доктора физико-математических наук Г.И. Синкевич [36]. Семинар Г.И. Синкевич в последнее время проводится дистанционно и привлекает много иностранных участников [37].

Имеется возможность принять в зарубежных семинарах по истории математики, проходящих на кафедре истории и философии науки в Кембриджском университете [38]. В настоящее время участники этих семинаров могут делать доклады (и слушать коллег) дистанционно, что делает их более доступными для исследователей всего мира.

Помимо годичной научной конференции ИИЕТ РАН им. Вавилова, проходящей в Москве и Санкт – Петербурге, на которой есть секции истории математики и механики, большое значение в последние годы приобрела международная конференция «Алгебра, теория чисел и дискретная геометрия: современные проблемы, приложения и проблемы истории», которая проходит ежегодно при Тульском государственном педагогическом университете. Кроме того, секция «история математики» с недавнего времени появилась на Крымской осенней математической школе (КРОМШ) и на Всероссийском симпозиуме по промышленной и прикладной математике (ВСППМ).

Следует отметить, что в нашей стране недавно было создано «Русское общество истории и философии науки» в целях усиления междисциплинарного взаимодействия в области науки, под эгидой которого проходят конференции и семинары [39].

Часть известных международных конференций по истории математики проходят в Лондоне и Париже, например «International Conference on History of Mathematics» [40]. Эти конференции призваны объединить ведущих ученых для обмена опытом и результатами исследований по всем аспектам истории математики. Они также обеспечивают междисциплинарную платформу для исследователей, практиков и преподавателей, чтобы представить и обсудить самые последние инновации, тенденции и проблемы, с которыми приходится сталкиваться в области истории математики.

Перейдём к теме подготовки докладов на международных конференциях.

Ряд учёных предпочитают делать доклады, используя только мел и доску. Однако, при слабом владении языком рекомендуется следует сделать презентацию в одной из версий Power Point. При этом желательно, чтобы текст был отредактирован специалистом. Во время показа следует проговаривать ту информацию, которая выводится на экран, пояснять ее. Вводную часть следует сделать короткой и сразу перейти к сути дела.

Если докладчик использует проектор, удобно пользоваться лазерной указкой. Это позволит привлечь внимание к ключевым моментам доклада. Если в докладе присутствуют громоздкие формулы, желательно использовать экран для их демонстрации. [41, с. 51-54]

После выступления на конференции (семинаре) с изложением своих результатов текст статьи, как правило, дополняется и претерпевает коррекцию в сторону улучшения.

*ж) Публикация.*

Заключительный этап исследований – представление своих результатов в виде, пригодном для опубликования. Так как разные журналы, предъявляют разные требования к электронной версии статьи, нужно сначала определиться с выбором журнала.

Рассмотрим возможные варианты:

1. Журнал «*История науки и техники*». Это ежемесячный научный журнал, который публикует научные материалы и обзоры в области истории науки (в том числе математики). Основным направлением журнала являются: публикация статей, в которых даётся анализ архивных документов и выявление ранее неизвестных фактов, представляющих историческую и научную ценность; публикация материалов, содержащих исторический анализ становления и развития науки, а также историю становления и развития научных школ и направлений и т.п. Индексируется в РИНЦ, входит в список ВАК. Текст статьи необходимо набирать в редакторе MS Word.
URL: http://int.tgizd.ru/.

2. Журнал «*Чебышевский сборник*». С 2015 г. принимает статьи по истории математики в формате TEX, в том числе большого объема. Среднее время нахождения статьи в редакции до момента публикации - 12 месяцев. С 2018 г. входит в базу данных Scopus и выходит четыре раза в год.
URL: https://www.chebsbornik.ru.

3. Журнал «*Научные ведомости БелГУ. Серия Прикладная математика & Физика*». Принимает статьи по истории математики в формате TEX. Среднее время нахождения статьи в редакции до момента публикации - шесть месяцев. Индексируется в РИНЦ, входит в список ВАК.
URL: http://nv.bsu.edu.ru/nv/mag/06/archive/.

4. Журнал «*Antiquitates Mathematicae*». Польский журнал, издаётся с 2007 г. Специализируется на истории и философии математики; принимает статьи в формате TEX на английском языке. Входит в список ВАК; использует перекрёстное рецензирование. URL: https://wydawnictwa.ptm.org.pl/index.php/antiquitates-mathematicae.

5. Журнал «Таврический вестник математики и информатики». Принимает статьи по истории математики в формате TEX. Индексируется в РИНЦ, входит в список ВАК. URL: https://tvim.info/

6. Журнал «*Известия высших учебных заведений. Прикладная нелинейная динамика*» - старейшее российское специализированное периодическое издание по нелинейной динамике, теории хаоса и их приложениям, издается с 1993 года. С 2008 года журнал входит в «Перечень периодических научных и научно-технических изданий РФ, рекомендуемых для публикации основных результатов диссертаций на соискание ученой степени доктора наук» ВАК РФ. С 2009 года индексируется в РИНЦ (загружены все выпуски, начиная с 1993 года). С 2017 года индексируется Scopus. С 2018 года индексируется Web of Science (ESCI). Принимает статьи в формате TEX.
URL: https://andjournal.sgu.ru/

7. Журнал «*Historia Mathematica*» - академический журнал по истории математики, издаваемый компанией Elsevier. Он был основан Кеннетом О. Май в 1971 г., как бесплатный информационный бюллетень Notae de Historia Mathematica, но к 1974 году превратился в полноценный журнал. Historia Mathematica публикует исследования по истории математики и ее развитию. В частности, журнал поощряет исследования о математиках и их работе в историческом контексте, об историографических темах в истории математики и о взаимосвязях между математическими идеями и наукой. Журнал индексируется в Scopus. Статьи принимаются на английском языке в формате TEX.
URL: https://www.journals.elsevier.com/historia-mathematica/

8. Журнал «*Archive for History of Exact Sciences*» принимает статьи высокого качества по истории математических наук. Его цель - обеспечить быструю и полную публикацию статей, а также пролить свет на концептуальную основу наук, анализируя исторический ход математической и естественнонаучной мысли. Журнал индексируется в Scopus. Статьи принимаются на английском языке.
URL: https://www.springer.com/journal/407

Так как значительная часть журналов принимает статьи в формате TEX, полезным будет ознакомление с издательской системой для набора текста LaTeX.

LaTeX – это система общего назначения для набора текстов, которая использует TEX в качестве средства форматирования. Благодаря своей гибкости, простоте использования и профессиональному полиграфическому качеству, LaTeX в настоящее время применяется при подготовке изданий почти по всем областям естественных наук. [42]

Система для набора текста LaTeX имеет несколько преимуществ перед текстовым редактором Word:

- Отличное качество подготовки научных текстов (удобно поддержана верстка математических формул).
- Автоматическая рубрикация документа, нумерация формул, рисунков и списка литературы. При этом все элементы грамотно располагаются на странице.
- TEX, форматирующее ядро LATEX, чрезвычайно мобилен и свободно доступен. Поэтому система работает практически на всех существующих платформах.

В освоении данный системы полезными будут пособия [43-44].

Многие ученые используют в своей работе веб-редактор LaTeX «Overleaf» . «Overleaf» — это бесплатная система совместного написания и опубликования текста, которая сильно ускоряет процесс создания научных работ как для авторов, так и издателей [45].

«Overleaf» можно охарактеризовать как сервис, который позволяет легко создавать, редактировать и делиться своими научными идеями в Интернете с помощью LaTeX.

### 2. Рекомендации по подготовке содержательной части научной статьи[4]

Для успешности статьи полезно рассмотреть ряд аспектов, которые помогут в ее написании. Желательно определить целевую аудиторию, от чего зависит глубина раскрытия основных и частных понятий. Она может быть адресована специалистам в данной области науки, математикам в целом и историкам математики[5].

Суть хорошей статьи по истории науки состоит в том, что она показывает независимую и критическую мысль. Статья не должна выглядеть, как отчет о какой-то теме, которая была затронута раньше: нужно найти такие способы представления материала, которые дают возможность для размышления.

Отправной точкой для развития мысли может послужить несогласие хотя бы с одним утверждением из вторичной литературы. Ошибки и неточности очень распространены, особенно в популярных книгах по истории математики. При проведении собственного исследования важно сосредоточиться на небольшом вопросе и попытаться выяснить, что различные вторичные источники говорят о нем. Как только вы хорошо поймете тему, скорее всего, обнаружится, что некоторые из более слабых вторичных источников очень поверхностны и, вполне возможно, совершенно неверны. Рекомендуется отметить такие недостатки в литературе и объяснить, что в них не так, и почему их ошибки значительны с точки зрения правильного понимания вопроса.

---

[4] Данный параграф написан с использованием материалов сайта доцента университета г. Утрехта (Нидерланды) В. Бласьо [46].
[5] Некоторые рекомендации П.Р. Халмоша, приведённые им в статье «Как писать математические тексты» [47], вполне подойдут и для подготовке статьи по истории математики.

Один из способов настроиться критически - это вступить в дискуссию, уже имеющуюся во вторичной литературе. Есть много случаев, когда историки науки расходятся во мнениях и предлагают конкурирующие интерпретации, часто в довольно жарких спорах. Выбор такой темы уведет от соблазна просто накапливать информацию и факты. Вместо этого придется критически взвесить доказательства и аргументы обеих сторон. Вероятно, вы окажетесь на той или иной стороне, и тогда будет вполне естественно внести свой собственный аргумент в пользу одной из сторон и свои собственные ответы на аргументы противоположной стороны.

Статья, построенная на сравнении и контрасте - это менее конфронтационный вариант статьи, основанной на дискуссии. Здесь мы также имеем дело с различными интерпретациями во вторичной литературе, но вместо того, чтобы пытаться "выбрать победителя", можно проиллюстрировать разнообразие подходов. Сравнивая различные точки зрения, мы поднимаем новые вопросы и освещаем новые стороны, которые не были очевидны, когда каждая точка зрения рассматривалась изолированно. Таким образом, более четко выявляются сильные и слабые стороны вопроса, а также предположения и следствия каждой точки зрения.

### 3. Получение академических и учёных степеней по истории математики

Поступление в магистратуру и аспирантуру дает возможность продолжить более глубокое изучение данной области и построить профессиональную карьеру в качестве исследователя истории математики.

Прием в аспирантуру по направлению «История науки и техники» осуществляется в Институте истории естествознания и техники им. С.И. Вавилова РАН (ИИЕТ РАН), город Москва. Для защиты диссертации по истории физико-математических наук в ИИЕТ РАН существует диссертационный совет, в состав которого входят доктора физико-математических, технических, исторических, философских и других наук. Соискателям, успешно прошедшим процедуру защиты диссертации, присваивается ученая степень кандидата (доктора) физико-математических наук по специальности «история науки и техники».

Для желающих продолжить своё образование в области истории и философии математики за границей можно порекомендовать следующие университеты:

- Кафедра истории и философии Университета Калгари, Канада
- Математический факультет Университета Саймона Фрейзера, Канада
- Институт истории и философии науки и техники Университета Торонто, Канада
- Исторический факультет Оксфордского университета, Великобритания
- Школа философии, религии и истории науки Лидского университета, Великобритания
- Кафедра истории науки Гарвардского университета, США
- Исторический факультет Чикагского университета, США
- Факультет философии Мичиганского университета, США.

### 4. Премии по истории математики

Самые успешные исследования в области истории математики, как это принято и в других областях науки, поощряются в виде различных премий и призов.

Рассмотрим наиболее известные из них.

**Taylor and Francis early career research prize (Премия Тейлора и Фрэнсиса за исследования в начале карьеры )**

Премия в размере 1000 фунтов стерлингов предоставляется за эссе объемом до 8000 слов по любому аспекту истории математики.

Победившее эссе принимается к публикации в журнале Британского Общества по истории математики. Победитель будет назначен координатором социальных сетей Общества в течение двух лет, получит доступ к аккаунтам Общества в Twitter и Facebook и будет приглашен вести блог по истории математики.

Претендентами на данную премию начинающие исследователи могут быть как аспиранты, так и лица, получившие степень кандидата наук не ранее, чем 5 лет назад [48].

**Конкурс диссертаций по истории науки DHST.**

Организаторами этого конкурса выступают Международный союз истории и философии науки и техники, Отдел истории науки и техники (IUHPST/DHST). Материалы диссертации, представляемой на конкурсе, должны быть посвящены истории науки, техники или медицины. Победителям предоставляется сертификат, освобождающий от регистрационных сборов конгрессе IUHPST/DHST и позволяющий компенсировать оплату проезда и проживания [49].

**Otto Neugebauer Prize (Премия Отто Нойгебауэра по истории математики)**

Премия присуждается за оригинальную и основополагающую работу в области истории математики, которая значительно углубляет наше понимание развития математики или конкретного математического предмета в любой период и в любом географическом регионе. Награда включает в себя именной сертификат и денежный приз в размере 5000 евро. Премия вручается на Европейском математическом конгрессе президентом Европейского математического общества. Лауреату дается возможность представить на конгрессе работу, получившую приз [50].

**Премия Неймана по истории математики**

Премия присуждается за монографию на английском языке (включая книги в переводе), посвященную истории математики и предназначенную для читателей-неспециалистов. Никаких дополнительных ограничений в отношении предмета, а также гражданства автора или места публикации не существует. Премия названа в честь Питера М. Неймана, бывшего президента и давнего члена Британского общества истории науки. Величина приза составляет 1000 фунтов стерлингов [51].

В 2009 году премия Неймана была присуждена Р. Нетц и У. Ноэль за работу «The Archimedes Codex: How a Medieval Prayer Book is Revealing the True Genius of Antiquity's Greatest Scientist» («Кодекс Архимеда: как средневековый молитвенник раскрывает истинного гения величайшего ученого Древности») [52].

**Kenneth O. May Prize in the History of Mathematics (Премия Кеннета О. Мэйя по истории математики).**

Премия Кеннета О. Мэйя вручается раз в 4 года за выдающийся вклад в историю математики. Впервые была присуждена в 1989 году Д. Дж. Струику и А.П. Юшкевичу (московскому учёному) на XVIII Международнov конгресса по истории науки. В 1993 году в рамках премии Кеннета О. Мэйя была учреждена бронзовая медаль [53].

**Премия Монтукла**

Премия Монтукла присуждается Исполнительным комитетом Международной комиссии по истории математики каждые четыре года молодому ученому - автору лучшей статьи, опубликованной в журнале Historia Mathematica за четыре года, предшествовавшие очередному Международному конгрессу по истории науки и техники. Премия составляет 1000 долларов [54]. Впервые премия Монтукла была присуждена в 2009 году.

**The Grattan-Guinness archival research travel grant (Грант на архивные исследования Grattan-Guinness)**

Гранты на поездки для проведения архивных исследований Grattan-Guinness были учреждены для оказания помощи ученым на ранних этапах их исследовательской карьеры в области истории и философии математики, а также в области истории математического образования и его связи с современными проблемами.

Указанные гранты открыты для докторантов или ученых, имеющих не более шести лет «постдоковских» исследований в области истории и/или философии математики. Гранты предоставляются специально для того, чтобы можно было добраться до места проведения исследований в архиве по выбору получателя. Такие гранты не предназначены для полного покрытия общей стоимости предлагаемого исследовательского проекта, но для покрытия командировочных расходов [55].

## Заключение

Авторы выражают надежду на то, что предложенная ими схема организации исследований в области истории математики и соответствующий инструментарий помогут молодым учёным

- правильно распределить свои силы на всех этапах исследований;
- структурировать свою научную работу от постановки задачи до отправки статьи в научный журнал;
- познакомиться с возможностями национальных библиотек России, Европы и США и использовать их для поиска первоисточников и статей историко-математического характера;
- получить представление об имеющихся всероссийских и международных конференциях и семинарах по истории математики;
- определиться с выбором журнала, в котором будет опубликована статья по истории математики;
- повысить уровень своей мотивированности за счёт возможного участия в конкурсах на соискание международной премий и грантов по истории науки;
- узнать, в какие учреждения можно поступить, для получения академической степени магистра или ученой степени кандидата наук по специальности «история математики и механики».

Ожидается, что в недалёком будущем часть рутинной работы, связанной с поиском и анализом первоисточников и материалов по близкой к исследуемой теме можно будет поручить электронному помощнику, работающему на принципах искусственного интеллекта [56], что позволит сэкономить время на содержательную часть исследований по истории науки.

Авторы выражают благодарность доценту В.П. Богатовой (Воронеж) за полезные обсуждения.

*Список литературы:*

1. Богатов Е.М., Коренев А.В., Михайлов И.С. Об организации научно-исследовательской работы по истории математики в современных условиях / Сб. научных статей 5-й Междунар. науч. конф. перспективных разработок молодых ученых «Наука молодых-будущее России» (10-11 декабря 2020 года), в 4-х томах, Том 1. Юго-Зап. гос. ун-т., - Курск:, 2020, с. 290-295.

2. Сайт онлайн-поддержки изучения курсов профессионально-исторической направленности для будущих и настоящих учителей математики. URL: http://pyrkov-professor.ru/ (дата обращения: 18.10.2020).

3. Математика в СССР за 30 лет (1917-1947). Курош А.Г. (гл. ред.), Маркушевич А.И. (ред.), Рашевский П.К. (ред.) // М.-Л.: ГИТТЛ, 1948. 1043 с.

4. Математика в СССР за 40 лет (1917-1957). Курош А.Г. (гл. ред.). Том 1. Обзорные статьи. Сборник М.: Физматгиз, 1959. 1002 с.

5. Математика в СССР за 40 лет (1917-1957). Курош А.Г. (гл. ред.). Том 2. Библиография. Сборник М.: Физматгиз, 1959. 821 с.

6. Лаврентьев М. А. Механика в СССР за 50 лет. В 4-х т. т. 1-3 // Л.И. Седов (гл. ред.), Зельдович Я.Б., Ишлинский А.Ю., Лаврентьев М.А., Михайлов Г.К., Мусхелишвили Н.И., Черный Г. Г. – 1968.

7. История отечественной математики. Отв. ред. Штокало И.З. В четырех томах. В пяти книгах. Том 1. - Киев: Наукова думка, 1966. - 493 с.

8. История отечественной математики. Отв. ред. Штокало И.З. В четырех томах. В пяти книгах. Том 2. - Киев: Наукова думка, 1967. - 616 с.

9. История отечественной математики. Отв. ред. Штокало И.З. В четырех томах. В пяти книгах. Том 3. - Киев: Наукова думка, 1968. - 726 с.

10. История отечественной математики. Отв. ред. Штокало И.З. В четырех томах. В пяти книгах. Том 4. Книга 1. Киев: Наукова думка, 1970. – 884 с

11. Jankvist U. T., Kjeldsen T. H. New avenues for history in mathematics education: Mathematical competencies and anchoring // Science & education. – 2011. – V. 20. – Iss. 9. – p. 831-862.

12. Сборник «Математическое Просвещение». Третья серия. URL: https://www.mccme.ru/free-books/matpros.html (дата обращения: 20.10.2020).

13. Тихомиров В.М. Геометрия выпуклости // Квант, 4 (2003). - с. 1-2.

14. Лемма Шпернера: приложения и обобщения. Занятие 1 URL: http://www.mathnet.ru/php/presentation.phtml?option_lang=rus&presentid=12048 (дата обращения: 16.03.2021).

15. Александров П. С. О некоторых основных направлениях в общей топологии // УМН, 19:6(120) (1964), с. 3–46.

16. Пырков В.Е. Тематический указатель статей сборника «Историко-математические исследования» за 1948 — 2009 годы / Ин-т истории естествознания и техники им. С.И. Вавилова РАН — М.: "Янус-К", 2011. — 84 с.

17. Gruber, P. M., & Wills, J. M. (Eds.). Handbook of convex geometry. North Holland. Elsever, 1993. 803 p.

18. MacTutor History of Mathematics archive. URL: https://mathshistory.st-andrews.ac.uk/ (дата обращения: 20.10.2020).

19. Борисов Ю. Ф., Залгаллер В. А. и др. К 90-летию со дня рождения А.Д. Александрова (1912–1999) // УМН, 57:5(347) (2002), с. 169–181

20. Канторович Л. В. Мой путь в науке // УМН, 1987. – Т. 42. – №. 2. – С. 183-213.

21. Стефан Банах (некролог) // УМН, 1:3-4(13-14) (1946), с.13–16

22. Архив Российской академии наук URL: http://www.arran.ru (дата обращения: 16.03.2021).

23. The Oxford Handbook of The History of Mathematics. Eleanor Robson& Jacqueline Stedall ed. Oxford: Oxford University Press, 2009. 918 p.

24. Grattan-Guinness I., ed . Companion Encyclopedia of the History and Philosophy of the Mathematical Sciences, 2 vols. Baltimore : J. Hopkins University Press, 2003. 1806 p.

25. Вилейтнер Г. История математики от Декарта до середины XIX столетия / Пер. с нем. под ред. А.П. Юшкевича. - М.: Физматлит, 1960. - 468 с.

26. Стройк Д. Я. Краткий очерк истории математики/ пер. с нем. и доп. И. Б. Погребысского. — 2-е изд. — М.: Наука, 1969. — 328 с.

27. James I. M. (ed.). History of topology. – North Holland. Amsterdam, 1999. - 1056 p.

28. Александров А.Д., Колмогоров А.Н., Лаврентьев М.А. (ред.) Математика: её содержание, методы и значение. Том 1. М.: Изд. АН СССР, 1956. – 296 с.
29. Александров А.Д., Колмогоров А.Н., Лаврентьев М.А. (ред.) Математика: её содержание, методы и значение. Том 2. М.: Изд. АН СССР, 1956. – 397 с.
30. Александров А.Д., Колмогоров А.Н., Лаврентьев М.А. (ред.) Математика: её содержание, методы и значение. Том 3. М.: Изд. АН СССР, 1956. – 336 с.
31. History of mathematics literature guide. URL: http://intellectualmathematics.com/blog/history-of-mathematics-literature-guide/ (дата обращения: 16.03.2021).
32. Электронная библиотека «Научное наследие России». URL: http://www.e-heritage.ru/about.html (дата обращения: 16.03.2021).
33. Project Euclid .URL: https://projecteuclid.org (дата обращения: 16.03.2021).
34. Lützen J. Joseph Liouville 1809–1882: Master of pure and applied mathematics. – Springer Science & Business Media, 2012. – V. 15, 885 p.
35. Историко-математические исследования / Под редакцией Г. Ф. Рыбкина и А. П. Юшкевича. — Выпуск XV. — М.: ГИФМЛ, 1963. — 480 с.
36. Одинец В.П., Синкевич Г.И. Сто докладов по истории математики на семинаре в Санкт-Петербурге // Математика в высшем образовании. 2018. - № 16. С. 79-84.
37. Семинары: Sumio Yamada, Evolution of geometric ideas, and the role of Relativity URL: http://www.mathnet.ru/php/seminars.phtml?option_lang=rus&presentid=30178 (дата обращения: 16.03.2021).
38. Department of History and Philosophy of Science URL: https://www.hps.cam.ac.uk/news-events/seminars-reading-groups/departmental (дата обращения: 16.03.2021).
39. РОИФН | Русское общество истории и философии науки. URL: http://rshps.org (дата обращения: 17.03.2021).
40. International Conference on History of Mathematics (ICHM) (URL: https://waset.org/history-of-mathematics-conference (дата обращения: 18.03.2021).
41. Сосинский А. Б. Mathematical English : Учебник английского для математиков. Издательство МЦНМО, 2017. 88 с.
42. Гуссенс М., Миттельбах Ф., Самарин А., Путеводитель по пакету Latex и его расширению Latex2e: Пер. с англ. – М.: Мир, 1999 – 621 с.
43. Балдин Е.М. Компьютерная типография LaTeX // БХВ-Петербург – 2012. 304 с.
44. Кнут Д. Все про TEX = The TEXBook. — М.: Вильямс, 2003 — 560 с.

45. Overleaf, Online LaTeX Editor URL: https://ru.overleaf.com/learn/latex/Russian (дата обращения: 18.03.2021).

46. How to write a history of mathematics essay URL: http://intellectualmathematics.com/blog/how-to-write-a-history-of-mathematics-essay/ (дата обращения: 19.03.2021).

47. Халмош П.Р. Как писать математические тексты // УМН, 1971, том 26, вып. 5(161), с. 243-269.

48. 2021 Taylor and Francis early career research prize URL: https://www.bshm.ac.uk/2021-taylor-and-francis-early-career-research-prize (дата обращения: 16.03.2021).

49. 2021 DHST dissertation prize competition call for applications URL: https://www.ichst2021.org/a2021-dhst-dissertation-prize-competition-call-for-applications/ (дата обращения: 16.03.2021).

50. Otto Neugebauer Prize URL: https://euro-math-soc.eu/otto-neugebauer-prize (дата обращения: 16.03.2021).

51. Neumann Prize URL: https://www.bshm.ac.uk/neumann-prize (дата обращения: 16.03.2021).

52. Mann T. History of Mathematics and History of Science // The University of Chicago Press on behalf of The History of Science Sosiety, 2011, p. 518-526.

53. Kenneth O. May Prize in the History of Mathematics URL: https://www.mathunion.org/ichm/prizes/kenneth-o-may-prize-history-mathematics (дата обращения: 16.03.2021).

54. ICHM Montucla Prize URL: https://www.mathunion.org/ichm/prizes/ichm-montucla-prize (дата обращения: 16.03.2021).

55. THE GRATTAN-GUINNESS ARCHIVAL RESEARCH TRAVEL GRANT URL: https://www.mathunion.org/fileadmin/ICHM/Grants/IGG-Grants-4.pdf (дата обращения: 16.03.2021).

56. Шамина О., Козлов Д. Автоматический поиск научных статей в сети Интернет // Материалы летней школы по информационному поиску RuSSIR 2008, Вып. 6. С. 43-62

**Bogatov Egor Mikhailovich,**

**Korenev Artem Viktorovich,**

**Mikhailov Ilya Andreevich**

**ABOUT THE MODERN TOOLS AND METHODS OF SCIENTIFIC RESEARCH CONDUCTING IN THE FIELD OF THE HISTORY OF MATHEMATICS**

*Abstract.* One of the variants for systematizing the activities of the historian of mathematics is proposed, as well as a scheme for organizing research and search work in the preparation of scientific articles and reports on the history of science.

*Keywords*: a scheme for organizing research on the history of mathematics, a methodology for preparing articles and reports on the history of mathematics, modern tools for historical and mathematical research; source base of the historian of mathematics.

**Bibliography**

1. Bogatov E.M., Korenev A.V., Mikhajlov I.S. About the organization of research work on the history of mathematics in modern conditions / Collection of scientific articles of the 5th International Scientific Conference of Promising Developments of Young Scientists "Young Science-the Future of Russia" (December 10-11, 2020), in 4 volumes, Vol. 1. Southwestern State University, Kursk:, 2020, p. 290-295.
2. Website of online support for the study of professional-historical courses for future and current teachers of mathematics. Available at: http://pyrkov-professor.ru/ (accessed: 18.10.2020).
3. Matematika v SSSR za 30 let (1917-1947). Ed. Kurosh A.G. // M.-L.: GITTL, 1948.
4. Matematika v SSSR za 40 let (1917-1957). Ed. Kurosh A.G. Book 1. Obzorny`e stat`i. Sbornik M.: Fizmatgiz, 1959. 1002 p.
5. Matematika v SSSR za 40 let (1917-1957). Ed. Kurosh A.G. Book 2. Bibliografiya. Sbornik M.: Fizmatgiz, 1959. 821 p.
6. Lavrentev M. A. Mekhanika v SSSR za 50 let. V 4-h t. T. 1-3 // L.I. Sedov (head editor), 1968.
7. Istoriya otechestvennoj matematiki. Ed. Shtokalo I.Z. In 4 volumes. In five books. Volume 1. - Kiev: Naukova dumka, 1966. - 493 p.
8. Istoriya otechestvennoj matematiki. Ed. Shtokalo I.Z. In 4 volumes. In five books. Volume 2. - Kiev: Naukova dumka, 1967. - 616 p.
9. Istoriya otechestvennoj matematiki. Ed. Shtokalo I.Z. In 4 volumes. In five books. Volume 3. - Kiev: Naukova dumka, 1968. - 726 p.
10. Istoriya otechestvennoj matematiki. Ed. Shtokalo I.Z. In 4 volumes. In five books. Volume 4. Book 1. Kiev: Naukova dumka, 1970. - 884 p.
11. Jankvist U. T., Kjeldsen T. H. New avenues for history in mathematics education: Mathematical competencies and anchoring //Science & education. – 2011. – V. 20. – Iss 9. – p. 831-862.
12. Matematicheskoe Prosveshchenie. Third series. Available at: https://www.mccme.ru/free-books/matpros.html (accessed: 20.10.2020).
13. Tikhomirov V. M. Geometry of convexity // Kvant, 4 (2003). - p.1-2.
14. Sperner's Lemma: applications and generalizations. Lesson 1. Available at: http://www.mathnet.ru/php/presentation.phtml?option_lang=rus&presentid=12048 (accessed: 16.03.2021).

15. Alexandrov P. S. On some basic directions in general topology // Uspekhi Mat. Nauk, 19:6(120) (1964), pp. 3–46.
16. Pyrkov V.E. «Istoriko-matematicheskie issledovaniya»: Tematicheskij ukazatel statej sbornika za 1948 — 2009 gody / V.E. Pyrkov; Ros. akad. nauk. In-t istorii estestvoznaniya i tekhniki im. S.I. Vavilova. — M.: "Yanus-K", 2011. — 84 p.
17. Gruber P. M., & Wills J. M. (Eds.). Handbook of convex geometry. North Holland. Elsever, 1993. 803 p.
18. MacTutor History of Mathematics archive.. Available at: https://mathshistory.st-andrews.ac.uk/ (accessed: 20.10.2020).
19. Borisov Y. F., Zalgaller V. A. and others, On the 90th anniversary of the birth of A. D. Alexandrov (1912–1999) // Uspekhi Mat. Nauk,, 57:5(347) (2002), pp. 169-181
20. Kantorovich L. V. My journey in science // Uspekhi mat. nauk. - 1987. - Vol. 42. - Iss 2. - p. 183-213.
21. Stefan Banach (obituary) // Uspekhi Mat. Nauk, 1:3-4(13-14) (1946), p. 13-16.
22. Archive of the Russian Academy of Sciences. Available at: http://www.arran.ru (accessed: 16.03.2021).
23. The Oxford Handbook of The History of Mathematics, edited by Eleanor Robson& Jacqueline Stedall. Oxford: Oxford University Press, 2009. 918 p.
24. Grattan-Guinness I., ed . Companion Encyclopedia of the History and Philosophy of the Mathematical Sciences, 2 vols. Baltimore : J. Hopkins University Press, 2003. 1806 p.
25. Wieleitner H. Geschichte der Mathematik, II Teil. Von Cartesius bis zur wende des 18. jahrhunderts, von dr. Heinrich Wieleitner: 1. hälfte. Arithmetik, algebra, analysis bearb. unter benutzung des nachlasses von dr. Anton von Braunmühl. 1911.-2. hälfte. Geometrie und trigonometrie. 1921; H. Wieleitner, Geschichte der Mathematik. II. Von 1700 bis zur Mitte des 19. Jahrhunderts. Berlin und Leipzig, 1923, p. 53-147.
26. Struick D. J. Abriss der Geschichte der Mathematik / Dirk J. Struick Berlin, 1963.
27. James I. M. (ed.). History of topology. – North Holland. Amsterdam, 1999. - 1056 p.
28. Alexandrov A.D., Kolmogorov A.N., Lavrentev M.A. (ed.) Matematika: eyo soderzhanie, metody i znachenie. Volume 1. M.: Izd. AN SSSR, 1956. – 296 p.
29. Alexandrov A.D., Kolmogorov A.N., Lavrentev M.A. (ed.) Matematika: eyo soderzhanie, metody i znachenie. Volume 2. M.: Izd. AN SSSR, 1956. – 397 p.
30. Alexandrov A.D., Kolmogorov A.N., Lavrentev M.A. (ed.) Matematika: eyo soderzhanie, metody i znachenie. Volume 3. M.: Izd. AN SSSR, 1956. – 336 p.
31. History of mathematics literature guide. Available at: http://intellectualmathematics.com/blog/history-of-mathematics-literature-guide/ (accessed: 16.03.2021).
32. Library «Nauchnoe nasledie Rossii». Available at: http://www.e-heritage.ru/about.html (accessed: 16.03.2021).
33. Project Euclid. Available at: https://projecteuclid.org (accessed: 16.03.2021).
34. Lützen J. Joseph Liouville 1809–1882: Master of pure and applied mathematics. Springer Science & Business Media, 2012. – V. 15, 885 p.
35. Istoriko-matematicheskie issledovaniya / Edited by G. F. Rybkin and A. P. Yushkevich. - Issue XV. - Moscow: GIFML, 1963. - 480 p.
36. Odinets V. P., Sinkevich G. I. One hundred reports on the history of mathematics at a seminar in St. Petersburg / / Mathematics in Higher Education. 2018. – Iss. 16. p. 79-84.
37. Seminars: Sumio Yamada, Evolution of geometric ideas, and the role of Relativity. Available at: http://www.mathnet.ru/php/seminars.phtml?option_lang=rus&presentid=30178 (accessed: 16.03.2021).

38. Department of History and Philosophy of Science. Available at: https://www.hps.cam.ac.uk/news-events/seminars-reading-groups/departmental (accessed: 16.03.2021).
39. RSHPS | Russian Society for the History and Philosophy of Science. Available at: http://rshps.org (accessed: 17.03.2021).
40. International Conference on History of Mathematics (ICHM). Available at: https://waset.org/history-of-mathematics-conference (accessed: 18.03.2021).
41. Sosinskii A. B. Mathematical English: An English textbook for mathematicians. Publishing house: MCCME, 2017, 87 p.
42. Goossens M., Mittelbach F. and Samarin A. The LaTeX companion. Vol. 1. Reading: Addison-Wesley, 1994, 528 p.
43. Baldin E. M. Computer printing house LaTeX // BHV-Petersburg-2012, 296 p.
44. Knuth D. E. The TEXbook Boston etc.: ADDISON-WESLEY PUBLISHING COMPANY, 1986, 483 p.
45. Overleaf, Online LaTeX Editor. Available at: https://ru.overleaf.com/learn/latex/Russian (accessed: 18.03.2021).
46. Halmosh P. R. How to write mathematical texts // Uspekhi Mat. Nauk, 1971, volume 26, issue 5 (161), p. 243-269.
47. How to write a history of mathematics essay. Available at: http://intellectualmathematics.com/blog/how-to-write-a-history-of-mathematics-essay/ (accessed: 19.03.2021).
48. 2021 Taylor and Francis early career research prize. Available at: https://www.bshm.ac.uk/2021-taylor-and-francis-early-career-research-prize (accessed: 16.03.2021).
49. 2021 DHST dissertation prize competition call for applications. Available at: https://www.ichst2021.org/a2021-dhst-dissertation-prize-competition-call-for-applications/ (accessed: 16.03.2021).
50. Otto Neugebauer Prize. Available at: https://euro-math-soc.eu/otto-neugebauer-prize (accessed: 16.03.2021).
51. Neumann Prize. Available at: https://www.bshm.ac.uk/neumann-prize (accessed: 16.03.2021).
52. Mann T. History of Mathematics and History of Science // The University of Chicago Press on behalf of The History of Science Society, 2011, p. 518-526.
53. Kenneth O. May Prize in the History of Mathematics. Available at: https://www.mathunion.org/ichm/prizes/kenneth-o-may-prize-history-mathematics (accessed: 16.03.2021).
54. ICHM Montucla Prize. Available at: https://www.mathunion.org/ichm/prizes/ichm-montucla-prize (accessed: 16.03.2021).
55. THE GRATTAN-GUINNESS ARCHIVAL RESEARCH TRAVEL GRANT. Available at: https://www.mathunion.org/fileadmin/ICHM/Grants/IGG-Grants-4.pdf (accessed: 16.03.2021).